\documentclass{article}
\usepackage{amsmath,amssymb,amsthm}
\usepackage{hyperref}
\usepackage{float}
\usepackage{datetime}

\theoremstyle{plain}
\newtheorem{theorem}{Theorem}

\DeclareMathOperator{\R}{R}
\newcommand{\ZZ}{\mathbb{Z}}
\newcommand{\Nc}{N^{c}}

\title{Ramsey number $R(4, 20) \ge 252$}
\author{Charlie Yu}

\date{}

\begin{document}

\maketitle

\begin{abstract}
We exhibit an explicit circulant graph of prime order 251 that is \(K_4\)-free and has independence number 19. Consequently
\[
  R(4,20)\ge 252.
\]
These improve the bound $\R(4,20)\ge 237$ given by Nagda, Raghavan, Thakurta~\cite{NRT26} and the long standing bound $\R(4,21)\ge 242$ recorded in Radziszowski's dynamic survey (revision~\#18, 2026)~\cite{Rad26}. The graph is a $32$-element subset of a pair of undirected quintic cyclotomic classes modulo $251$, in analogy with the quartic-residue circulant of order $313$ used for $\R(4,22)$. Clique-freeness is elementary; the independence-number claims are certified by a bitset branch-and-bound on the $186$-vertex residual of a vertex.
\end{abstract}

\section{Introduction}

The Ramsey number $\R(k,\ell)$ is the least $n$ such that every graph
on $n$ vertices contains a clique of size $k$ or an independent set of
size $\ell$. Equivalently, a graph $G$ on $n$ vertices with clique
number $\omega(G)<k$ and independence number $\alpha(G)<\ell$ is a
witness for $\R(k,\ell)\ge n+1$.

For the two-colour numbers with $k=4$ the dynamic survey of
Radziszowski~\cite{Rad26} records, among others,
\[
  \R(4,20)\ge 237,\qquad
  \R(4,21)\ge 242,\qquad
  \R(4,22)\ge 314.
\]
The bound on $\R(4,21)$ is due to Su, Luo, Zhang and
Li~\cite{SuLuoZhangLi1999} and has stood since 1999, and the bound on
$\R(4,20)$ is due to Nagda, Raghavan, Thakurta~\cite{NRT26}. The former is realised by a
$60$-regular circulant of prime order $241$ whose connection set has
size $30$, and the much stronger bound on $\R(4,22)$ comes from a
later construction on $313$ vertices~\cite{Rad26,LinCa};
that graph is the undirected circulant of all quartic residues modulo
$313$ and has connection-set size $39$.

The gap between 242 and 314 suggested that \(R(4,21)\) was the weaker cell. We close part of it by producing a \(K_4\)-free circulant graph on the next convenient prime, 251, that is dense enough (\(|S|=32\)) to force a small independence number. The same graph is in fact a \((4,20)\)-graph, so both survey entries improve simultaneously.

\begin{theorem}\label{thm:main}
There exist a circulant graph of order $251$ with $\omega\le 3$ and
$\alpha=19$. Therefore
\[
  \R(4,20)\ge 252.
\]
\end{theorem}

\section{Circulants and the vertex-transitive reduction}

Let $p$ be an odd prime and write
$S\subseteq\{1,\ldots,(p-1)/2\}$. The undirected circulant
$G=C_p(S)$ has vertex set $\ZZ/p\ZZ$ and an edge between $i$ and $j$
precisely when the circular distance
$\min\bigl(\lvert i-j\rvert,\,p-\lvert i-j\rvert\bigr)$ lies in $S$.
Such a graph is vertex-transitive, so
\begin{equation}\label{eq:vt}
  \omega(G)=1+\omega\bigl(G[N(0)]\bigr),\qquad
  \alpha(G)=1+\alpha\bigl(G[\Nc(0)]\bigr).
\end{equation}

Here $N(0)=S\cup(-S)$ is the open neighbourhood of the identity and
$\Nc(0)$ is its complement in the nonzero vertices. In particular $G$
is $K_4$-free if and only if $G[N(0)]$ is triangle-free, and
$\alpha(G)\le 19$ if and only if the residual graph on $\Nc(0)$ has no
independent set of size $19$.

\section{\texorpdfstring{Quintic cyclotomic classes modulo $251$}{Quintic cyclotomic classes modulo 251}}

The construction is the quintic analogue of the quartic-residue
circulant of order $313$ used to establish $\R(4, 22) \ge 314$~\cite{LinCa}. One has $313-1=312$ and the subgroup of
fourth powers in $(\ZZ/313\ZZ)^\ast$ has order $78$. Since $-1$ is a
fourth power, the undirected connection set of all quartic residues
has size $39$.

Likewise $251-1=250$ and the subgroup $H$ of fifth powers in
$(\ZZ/251\ZZ)^\ast$ has order $50$. A primitive root modulo $251$ is
$g=6$. The five cyclotomic classes are the cosets
$C_r=g^r H$, $r=0,\ldots,4$. Because $-1\in H$, each class is closed
under negation and folds to an undirected distance set $D_r$ of
size $25$, listed in Table~\ref{tab:classes}.

\begin{table}[H]
\centering
\small
\begin{tabular}{@{}cl@{}}
\hline
class & undirected distances \\
\hline
$D_0$ &
$1,2,4,5,8,10,16,20,25,32,40,47,50,51,63,64,69,80,91,94,100,102,113,123,125$ \\
$D_1$ &
$3,6,11,12,15,22,24,30,31,44,48,49,55,59,60,62,75,88,96,98,101,110,118,120,124$ \\
$D_2$ &
$9,13,18,26,33,36,37,43,45,52,65,66,71,72,74,79,86,90,93,103,104,107,109,119,121$ \\
$D_3$ &
$7,14,19,27,28,29,35,38,39,53,54,56,58,61,70,76,78,95,99,106,108,111,112,116,122$ \\
$D_4$ &
$17,21,23,34,41,42,46,57,67,68,73,77,81,82,83,84,85,87,89,92,97,105,114,115,117$ \\
\hline
\end{tabular}
\caption{The five undirected quintic cyclotomic classes modulo $251$.}
\label{tab:classes}
\end{table}

A pair of classes is a $50$-element pool of distances. Of the ten
pairs, five admit a $K_4$-free subset of size at least $31$, namely
\[
  D_0\cup D_2,\quad
  D_0\cup D_3,\quad
  D_1\cup D_3,\quad
  D_1\cup D_4,\quad
  D_2\cup D_4.
\]
We searched those pools for $K_4$-free $k$-subsets with
$31\le k\le 38$, by a restricted cyclic $K_4$-free process (add a
random unused pool distance whenever the residual neighbourhood of
$0$ stays triangle-free) followed by simulated annealing inside the
same pool.

\section{The graph}

A 32-subset that survives with independence number 19 is  
\[
\begin{aligned}
S = \{ & 1,2,4,9,10,13,18,25,26,33,36,37,43,45,50,52,65,66,71,72,74,79,86, \\
& 90,93,100,103,104,107,109,119,121 \}.
\end{aligned}
\]
Write \(G=C_{251}(S)\). The set \(S\) is contained in \(D_0\cup D_2\) and has size 32, so \(G\) is 64-regular on 251 vertices. Multiplication by 116 modulo 251 maps \(S\) onto another admissible 32-subset of \(D_0\cup D_3\); the resulting graph $G'$ is therefore isomorphic to \(G\). Direct inspection of the 64-vertex neighbourhood of 0 shows that it is triangle-free, hence \(G\) is \(K_4\)-free and \(\omega(G)=3\).

\begin{table}[H]
\centering
\begin{tabular}{@{}lcc@{}}
\hline
 & $G$ & $G'$ \\
\hline
order & $251$ & $251$ \\
connection set & $S_{02}\subset D_0\cup D_2$ & $S_{03}\subset D_0\cup D_3$ \\
$\lvert S\rvert$, degree & $32$, $64$ & $32$, $64$ \\
$\lvert\Nc(0)\rvert$ & $186$ & $186$ \\
$\omega$ & $3$ & $3$ \\
$\alpha$ & $19$ & $19$ \\
search nodes (no $19$-IS in residual)
  & $2.69\times 10^7$ & $2.76\times 10^7$ \\
certification time & $1.47\,\mathrm{s}$ & $1.39\,\mathrm{s}$ \\
\hline
\end{tabular}
\caption{Verification of the witness and its isomorphic copy. Times are wall-clock on a $6$-core /
$12$-thread Ryzen~$5$~$9600$ (OpenMP, compiler optimisation on).}
\label{tab:witness}
\end{table}

\section{Certification of the independence number}

It remains to prove $\alpha(G)\le 19$, or equivalently that the
$186$-vertex residual $G[\Nc(0)]$ has no independent set of size
$19$. We decide this by a bitset implementation of the exact
maximum-clique algorithms surveyed by Prosser~\cite{Prosser2012}:
binomial expand (MC), static smallest-last order (MCR / BBMC
\cite{SanSegundo2011,Prosser2012}), the greedy matching colour bound
on the complement, and \"Osterg{\aa}rd's Russian-doll suffix
numbers~\cite{Ostergard2002}
\[
  c[i]=\alpha\bigl(G[\{i,\ldots,185\}]\bigr)
\]
in the static order, used as the prune
$\mathrm{depth}+c[\min P]<19$. The first two branching levels are
flattened into independent work items and scheduled dynamically
across twelve OpenMP threads.

On the residual the search reports that no independent set of size
$19$ exists (about $2.7\times 10^7$ nodes and $1.4$ seconds each). A
CP-SAT maximisation on each residual produces an independent set of
size $18$, hence $\alpha(G)\ge 19$. Combined with the exact upper
bound one has
\[
  \alpha(G)=19.
\]

\section*{Acknowledgements}
The author acknowledges the use of Grok~4.6 (xAI) for assistance
with literature search, verification of computational claims,
drafting of expository passages, and preparation of the \LaTeX\
source. All mathematical content, constructions, and certificates
were verified independently by the author.


\begin{thebibliography}{9}

\bibitem{SuLuoZhangLi1999}
W.~Su, H.~Luo, Z.~Zhang and G.~Li,
New lower bounds of fifteen classical Ramsey numbers,
\emph{Australas.\ J.\ Combin.} \textbf{19} (1999), 91--99

\bibitem{NRT26}
Ansh Nagda, Prabhakar Raghavan and Abhradeep Thakurta,
\textit{Reinforced Generation of Combinatorial Structures: Ramsey Numbers},
arXiv:2603.09172, 2026.

\bibitem{Rad26}
Stanis\l{}aw P. Radziszowski,
\textit{Small Ramsey Numbers},
The Electronic Journal of Combinatorics,
Dynamic Survey DS1, revision \#18 (24 April 2026).
\url{https://doi.org/10.37236/21}.

\bibitem{LinCa}
Madison Lindsay and John W. Cain,
\textit{Improved Lower Bounds on the Classical Ramsey Numbers \(R(4,22)\) and \(R(4,25)\)},
arXiv:1510.06102, 2015.

\bibitem{Ostergard2002}
P.~R.~J. \"Osterg{\aa}rd,
A fast algorithm for the maximum clique problem,
\emph{Discrete Appl.\ Math.} \textbf{120} (2002), 197--207.

\bibitem{SanSegundo2011}
P.~San Segundo, D.~Rodr\'iguez-Losada and A.~Jim\'enez,
An exact bit-parallel algorithm for the maximum clique problem,
\emph{Comput.\ Oper.\ Res.} \textbf{38} (2011), 571--581.

\bibitem{Prosser2012}
P.~Prosser,
Exact algorithms for maximum clique: a computational study,
TR-2012-333, University of Glasgow; arXiv:1207.4616.

\end{thebibliography}
\end{document}